\documentclass[11pt,a4paper]{article}
\usepackage{pifont}
\usepackage{bbding}
\usepackage{amssymb}
\usepackage{mathrsfs}
\usepackage{amsfonts}
\usepackage{amsthm,amscd,amsmath}
\usepackage[all]{xy}
\usepackage{pmat}
\usepackage{arydshln}
\usepackage{enumitem}

\DeclareMathOperator{\sgn}{{\rm sgn}}

\makeatletter
\def\@biblabel#1{#1}
\makeatother

\newtheorem*{theorem1}{Theorem}
\newtheorem{theorem}{Theorem}[section]
\newtheorem{lemma}[theorem]{Lemma}

\theoremstyle{definition}
\newtheorem{definition}[theorem]{Definition}

\theoremstyle{proposition}
\newtheorem{proposition}[theorem]{Proposition}

\theoremstyle{remark}
\newtheorem{remark}[theorem]{Remark}

\theoremstyle{corollary}
\newtheorem{corollary}[theorem]{Corollary}
\date{}

\begin{document}

\title{A New Proof of Sturm's Theorem via Matrix Theory}

\author{Kaiwen Hou \\
\small{Tsinghua University, P.R.China} \\
\small{hkw18@mails.tsinghua.edu.cn}
\and Bin Li \\
\small{Wuhan Foreign Languages School, P.R.China} \\
\small{li\_bin0529@sina.com}
}

\maketitle

\begin{abstract}
 By the classical Sturm's theorem, the number of distinct real roots of a given real polynomial $f(x)$ within any interval $(a,b]$ can be expressed by the number of variations in the sign of the Sturm chain at the bounds. Through constructing the ``Sturm matrix", a symmetric matrix  associated with $f(x)$ over $\mathbb R[x]$, variations in the sign of $f(x)$ can be characterized by the negative index of inertia. Therefore, this paper offers a new proof of Sturm's theorem using matrix theory.
\end{abstract}


\section{Introduction}

By the fundamental theorem of algebra, every non-constant polynomial $f(x)$ over the complex field $\mathbb C$ has a root in $\mathbb C$, furthermore, the number of complex roots counted with multiplicities equals the degree of $f(x)$.
Correspondingly, as a powerful tool to study the real roots of a given polynomial $f(x)\in\mathbb R[x]$, the classical Sturm's theorem not only yields the number of distinct real roots of $f(x)$, but also locates them in intervals.

For a non-constant polynomial $f(x)\in\mathbb R[x]$, let $f_0(x)=f(x)$ and $f_1(x)=f^{\prime}(x)$, the derivative of $f(x)$. The \textbf{canonical Sturm chain} $f_0(x), f_1(x), \ldots,$\\
$f_m(x)$ is obtained by the following modified Euclidean algorithm,
\begin{equation*}
\begin{split}
f_0(x) &= d_1(x)f_1(x)-f_2(x),\quad  \deg(f_2)<\deg(f_1),\\
f_1(x) &= d_2(x)f_2(x)-f_3(x),\quad  \deg(f_3)<\deg(f_2),\\
&\;\; \vdots\\
f_{m-2}(x) &= d_{m-1}(x)f_{m-1}(x)-f_m(x),\quad  \deg(f_{m})<\deg(f_{m-1}),\\
f_{m-1}(x) &= d_{m}(x)f_{m}(x).
\end{split}
\end{equation*}
For a sequence of non-zero real numbers $c_1, c_2, \ldots, c_n$, we denote by $\sigma(c_1,\ldots, c_n)$ the number of variations in sign of the sequence, i.e.
\begin{equation*}
\sigma(c_1,\ldots, c_n)=\sum_{i=1}^{n-1}\delta_{-1,\sgn(c_ic_{i+1})},
\end{equation*}
where $\delta$ is the Kronecker delta and $\sgn$ denotes the sign function. For a general sequence of real
numbers $c_1, \ldots, c_n$, let $c_{i_1},\ldots, c_{i_m}$ be the subsequence of it obtained by deleting all the zeros and define
 \begin{equation*}
\sigma(c_1,\ldots, c_n)=\sigma(c_{i_1},\ldots, c_{i_m}).
\end{equation*}
Set $V_f(x)=\sigma(f_0(x),f_1(x),\ldots,f_m(x))$.
\begin{theorem1}[Sturm's Theorem\cite{jacobson}]
For a non-constant polynomial $f(x)\in \mathbb R[x]$ and two real numbers $a<b$, the number of distinct real roots in the interval $(a,b]$ is $V_f(a)-V_f(b)$ if neither $a$ nor $b$ is a multiple root of $f(x)$.
\end{theorem1}
A matrix explanation of $V_f(x)$ will be given as follows and thus a new proof of the theorem above is obtained.

\section{Preliminaries}
Given a real symmetric matrix $A$, the quadratic form defined by $A$ can be brought to a diagonal form by a non-singular transformation of coordinates. The famous Sylvester's law of inertia states that the number of positive (negative)  coefficients in the diagonal form is an invariant which is called the positive (negative) index of inertia of $A$, denoted by $p(A)$ (respectively, $q(A)$). Note that $p(A)$ and $q(A)$ are the numbers of positive and negative eigenvalues (counted with multiplicities) of $A$ respectively.

\begin{lemma}~\label{lem1}
If $A$ is an $n\times n$ real symmetric matrix of rank $r(A)\geqslant n-1$ and  $B=\left(
                                                                                \begin{array}{cc}
                                                                                  A & \alpha \\
                                                                                  \alpha^T & b \\
                                                                                \end{array}
                                                                              \right)
$ with $\alpha\in\mathbb R^n$ and $b\in\mathbb R$ such that their determinants $|A|$ and $|B|$ satisfy $|A|^2+|B|^2\neq 0$, then
\begin{displaymath}
q(B)= \left\{ \begin{array}{ll}
q(A) & \textrm{if $|A||B|>0$ or $|B|=0$}\\
q(A)+1 & \textrm{if  $|A||B|<0$ or $|A|=0$.}
\end{array} \right.
\end{displaymath}
\end{lemma}
\begin{proof}
It is easy to prove the case when $|A||B|\neq 0$. If $|B|=0$, the condition $|A|^2+|B|^2\neq 0$ implies that $A$ is invertible.
We have
\begin{equation*}
\left(
                                                                                \begin{array}{cc}
                                                                                  I & 0 \\
                                                                                  -\alpha^TA^{-1} & 1 \\
                                                                                \end{array}
                                                                              \right)\left(
                                                                                \begin{array}{cc}
                                                                                  A & \alpha \\
                                                                                  \alpha^T & b \\
                                                                                \end{array}
                                                                              \right)\left(
                                                                                \begin{array}{cc}
                                                                                  I & -A^{-1}\alpha \\
                                                                                  0 & 1 \\
                                                                                \end{array}
                                                                              \right)=\left(
                                                                                \begin{array}{cc}
                                                                                  A & 0 \\
                                                                                  0 & b-\alpha^TA^{-1}\alpha\\
                                                                                \end{array}
                                                                              \right).
\end{equation*}
By taking the determinant on both sides, we obtain that $b-\alpha^TA^{-1}\alpha=0$. Thus $B$ conjugates to
\begin{equation*}
\left(
                                                                                \begin{array}{cc}
                                                                                  A & 0 \\
                                                                                  0 & 0 \\
                                                                                \end{array}
                                                                              \right),
\end{equation*}
 which implies $q(B)=q(A)$. If $|A|=0$, then $r(A)=n-1$ and $|B|\neq 0$. There exists an orthogonal matrix $Q$ such that
 \begin{equation*}
 Q^TAQ=\left(
                                                                                \begin{array}{cccc}
                                                                                  \lambda_1&  & & \\
                                                                                           &\ddots& &  \\
                                                                                           &      &\lambda_{n-1}& \\
                                                                                           &      &             &0\\
                                                                                \end{array}
                                                                              \right)
 \end{equation*}
 in which $\lambda_i$ ($1\leqslant i\leqslant n-1$) are all the nonzero eigenvalues of $A$. One has
\begin{equation*}
\left(
                                                                                \begin{array}{cc}
                                                                                  Q^T &  \\
                                                                                   & 1 \\
                                                                                \end{array}
                                                                              \right)B\left(
                                                                                \begin{array}{cc}
                                                                                  Q &  \\
                                                                                   & 1 \\
                                                                                \end{array}
                                                                              \right)
                                                                              =\begin{pmat}({...|})
\lambda_1 &        &               &   &b_1\cr
          & \ddots &               &   &b_2\cr
          &        & \lambda_{n-1} &   &\vdots\cr
          &        &               & 0 &b_n\cr\-
b_1       &b_2     &\cdots         &b_n&b\cr
\end{pmat}
\end{equation*}
which conjugates to a matrix of the form
\begin{equation*}
\begin{pmat}({...|})
\lambda_1 &        &               &   & \cr
          & \ddots &               &   &  \cr
          &        & \lambda_{n-1} &   &  \cr
          &        &               & 0 &c\cr\-
          &        &               & c &d\cr
\end{pmat}.
\end{equation*}
Since $|B|\neq 0$, $c\neq 0$ and thus the two eigenvalues $\mu_1, \mu_2$ of $\left(
                                                                                \begin{array}{cc}
                                                                                  0 & c \\
                                                                                  c & d \\
                                                                                \end{array}
                                                                              \right)$ satisfy
\begin{equation*}
\mu_1\mu_2=-c^2<0.
\end{equation*}
This implies $q(B)=q(A)+1$.
\end{proof}

Given an $n\times n$ matrix $A$ and $I\subseteq \{1,\ldots, n\}$, let $A_I$ denote the principal submatrix of $A$ whose rows and columns are indexed by $I$.
The determinant of $A_I$ is called the principal minor of $A$ of order $\# I$, with $\# I$ denoting the cardinality of $I$.
\begin{definition} Let $A$ be an $n\times n$ matrix over $\mathbb R$,
\begin{itemize}
\item[(i)] we call $D_1, D_2, \ldots, D_n$ a \textbf{principal minor sequence} of $A$, if $D_i=|A_{I_i}|$ such that $I_i\subset I_{i+1}\subseteq\{1, \ldots, n\}$ and $\#I_i=i$;
\item[(ii)] A principal minor sequence $D_1, D_2, \ldots, D_n$ is said to be \textbf{normal} if $D_r\neq 0$ with $r=r(A)$ and for $1\leqslant i\leqslant r-1$, any two consecutive minors $D_i, D_{i+1}$ are not both zero.
\end{itemize}
\end{definition}
Applying Lemma~\ref{lem1}, we obtain the following lemma.
\begin{lemma}~\label{lem2}
Given an $n\times n$ real symmetric matrix $A$ of rank $r=r(A)$, if $D_1, \ldots, D_n$ is a normal principal minor sequence of $A$, then
\begin{equation*}
q(A)=\delta_{-1,\sgn(D_1)}+\sum_{i=1}^{r-1}\delta_{-1,\sgn(D_iD_{i+1})}+\sum_{i=1}^{r-1}\delta_{0,D_i}.
\end{equation*}
\end{lemma}
\begin{remark}
It is worth pointing out here that one can prove, for any real symmetric matrix $A$, there exists a normal principal minor sequence by which we can determine its negative as well as positive indices of inertia.
\end{remark}

\section{Sturm matrix and Sturm's theorem}

We assume throughout that $(f(x),g(x))$ is a pair of non-constant real polynomials unless otherwise stated.

\subsection{Sturm matrix and its negative index of inertia }
Let $f_0(x)=f(x)$ and $f_1(x)=g(x)$. The \textbf{Sturm chain} $f_0, \ldots, f_m$ associated to $(f,g)$ is obtained as follows
\begin{equation*}
\begin{split}
f_0(x) &= d_1(x)f_1(x)-f_2(x),\quad  \deg(f_2)<\deg(f_1),\\
f_1(x) &= d_2(x)f_2(x)-f_3(x),\quad  \deg(f_3)<\deg(f_2),\\
&\;\; \vdots\\
f_{m-2}(x) &= d_{m-1}(x)f_{m-1}(x)-f_m(x),\quad  \deg(f_{m})<\deg(f_{m-1}),\\
f_{m-1}(x) &= d_{m}(x)f_{m}(x).
\end{split}
\end{equation*}
The following facts about the Sturm chain are obvious.
\begin{itemize}
\item[(a)] If $a$ is a common root of $f(x)$ and $g(x)$, then $f_i(a)=0$ for all $0\leqslant i\leqslant m$; otherwise, any two consecutive polynomials $f_i, f_{i+1}$ in the chain do not have a common root.
\item[(b)] $f_{m}(x)=gcd(f(x),g(x))$, hence $a$ is a common root of $f(x)$ and $g(x)$ if and only if $f_{m}(a)=0$.
\item[(c)] If $f_i(a)=0$ for some $1\leqslant i\leqslant m-1$ and $f_m(a)\neq 0$, then $f_{i-1}(a)f_{i+1}(a)<0$.
\end{itemize}
\begin{definition}
For $f(x), g(x)\in\mathbb R[x]$, the \textbf{Sturm matrix} associated to $(f(x),g(x))$ is defined to be
\begin{equation*}
\left(
                                                                                \begin{array}{cccc}
                                                                                  d_1(x)& 1 & & \\
                                                                                  1     &d_2(x)&\ddots &  \\
                                                                                           & \ddots    &\ddots&1 \\
                                                                                           &      &    1    &d_m(x)\\
                                                                                \end{array}
                                                                              \right),
\end{equation*}
denoted by $S_{f,g}(x)$.
\end{definition}
Note that $r(S_{f,g}(x))\geqslant m-1$ and for any common divisor $d(x)$ of $f(x)$ and $g(x)$,
\begin{equation*}
S_{f,g}(x)=S_{\frac{f}{d},\frac{g}{d}}(x).
\end{equation*}
We define the \textbf{refined Sturm chain} associated to $(f(x),g(x))$ by
\begin{equation*}
\tilde{f}_0, \tilde{f}_1, \ldots,\tilde{f}_m,
\end{equation*}
where $\tilde{f}_i=\frac{f_i}{f_m}$.  Especially we have
$S_{f_0,f_1}(x)=S_{\tilde{f}_0,\tilde{f}_1}(x)$.
Let $V_{f,g}(x)$ denote the number of variations in sign of $f_0(x), \ldots, f_m(x)$, i.e.
\begin{equation*}
V_{f,g}(x)=\sigma(f_0(x), \ldots, f_m(x)).
\end{equation*}
See that
\begin{equation*}
V_{f,g}(a)=V_{\tilde{f},\tilde{g}}(a)
\end{equation*}
with $\tilde{g}=\tilde{f}_1$ if $a$ is not a common root of $f(x)$ and $g(x)$.

We consider the following principal minor of $S_{f,g}(x)$,
\begin{equation*}
D_i(x)=|S_{f,g}(x)_{\{m-i+1,m-i+2,\ldots, m\}}|,
\end{equation*}
which is of order $i$ for $1\leqslant i\leqslant m$.
\begin{proposition}~\label{minor}
For $1\leqslant i\leqslant m$, $D_i(x)=\tilde{f}_{m-i}$.
\end{proposition}
\begin{proof} We use induction on $i$. See that
$D_1(x)=d_m(x)=\frac{f_{m-1}(x)}{f_m(x)}=\tilde{f}_{m-1}(x)$ and
\begin{equation*}
\begin{split}
D_2(x)&=\left|
                                                                                \begin{array}{cc}
                                                                                  d_{m-1}(x)& 1 \\
                                                                                  1     &d_m(x)\\
                                                                                \end{array}
                                                                              \right|=d_{m-1}(x)d_m(x)-1\\
     &=\frac{f_{m-1}(x)d_{m-1}(x)}{f_m(x)}-1=\frac{f_{m-2}(x)+f_m(x)}{f_m(x)}-1\\
     &=\frac{f_{m-2}(x)}{f_m(x)}=\tilde{f}_{m-2}(x).
\end{split}                                                                              \end{equation*}
Assume that $D_k(x)=\tilde{f}_{m-k}$ and $D_{k+1}(x)=\tilde{f}_{m-k-1}$.
\begin{equation*}
D_{k+2}(x)=\left|
                                                                                \begin{array}{cccc}
                                                                                  d_{m-k-1}(x)& 1 & & \\
                                                                                  1     &d_{m-k}(x)&\ddots &  \\
                                                                                           & \ddots    &\ddots&1 \\
                                                                                           &      &    1    &d_m(x)\\
                                                                                \end{array}
                                                                              \right|.
\end{equation*}
Expanding the determinant by the first row, we have
\begin{equation*}
\begin{split}
D_{k+2}(x)&=d_{m-k-1}(x)D_{k+1}(x)-D_k(x)=d_{m-k-1}(x)\tilde{f}_{m-k-1}(x)-\tilde{f}_{m-k}(x)\\
          &=\frac{d_{m-k-1}(x)f_{m-k-1}(x)-f_{m-k}(x)}{f_m(x)}=\tilde{f}_{m-k-2}(x)
\end{split}
\end{equation*}
which completes the proof.
\end{proof}
\begin{theorem}~\label{var}
$q(S_{f,g}(a))=V_{f,g}(a)$ if $a$ is not a common root of $f(x)$ and $g(x)$.
\end{theorem}
\begin{proof}
By the facts (a), (b) about the (refined) Sturm chain, $D_1(a), \ldots, D_m(a)$ is a normal principal minor sequence of $S_{f,g}(a)$ if $a$ is not a common root of $f(x)$ and $g(x)$. Hence by Lemma~\ref{lem2}, one has
\begin{equation}~\label{q1}
q(S_{f,g}(a))=\delta_{-1,\sgn(D_1(a))}+\sum_{i=1}^{r-1}\delta_{-1,\sgn(D_i(a)D_{i+1}(a))}+\sum_{i=1}^{r-1}\delta_{0,D_i(a)}.
\end{equation}
with $r=r(S_{f,g}(a))$.
Notice that $r(S_{f,g}(a))=m$ if $a$ is not a root of $f(x)$, otherwise, $r(S_{f,g}(a))=m-1$. In both cases, (\ref{q1}) can be rewritten as
\begin{equation}
q(S_{f,g}(a))=\delta_{-1,\sgn(D_1(a))}+\sum_{i=1}^{m-1}\delta_{-1,\sgn(D_i(a)D_{i+1}(a))}+\sum_{i=1}^{m-1}\delta_{0,D_i(a)}.
\end{equation}
Thus
\begin{equation*}
\begin{split}
q(S_{f,g}(a))&=\delta_{-1,\sgn(\tilde{f}_{m-1}(a))}+\sum_{i=1}^{m-1}\delta_{-1,\sgn(\tilde{f}_{m-i}(a)\tilde{f}_{m-i-1}(a))}+\sum_{i=1}^{m-1}\delta_{0,\tilde{f}_{m-i}(a)}\\
&=\sigma(1, \tilde{f}_{m-1}(a), \tilde{f}_{m-2}(a),\ldots,\tilde{f}_0(a))\\
&=\sigma(\tilde{f}_0(a), \ldots, \tilde{f}_{m}(a))\\
&=V_{f,g}(a),\\
\end{split}
\end{equation*}
in which the second equality follows from the fact (c) about the Sturm chain.
\end{proof}
\subsection{Proof of Sturm's theorem}
Let $a_1, \ldots, a_k\in\mathbb R$ be all the roots of $\tilde{f}_0(x)$. Suppose that
\begin{equation*}
a_1<a_2<\ldots<a_k.
\end{equation*}
Let $I_i$ be the open interval $(a_i, a_{i+1})$ for $0\leqslant i\leqslant k$, with the convention $a_0=-\infty$ and $a_{k+1}=+\infty$.
\begin{proposition}~\label{constant}
$q(S_{f,g}(x))$ is constant in each interval $I_i$.
\end{proposition}
\begin{proof}
Since $S_{f,g}(x)$ is a real symmetric matrix for every given $x\in\mathbb R$, the eigenvalues of $S_{f,g}(x)$
\begin{equation*}
\lambda_1(x), \ldots, \lambda_m(x)
\end{equation*}
are all real numbers.
Suppose $\lambda_1(x)\leqslant \ldots\leqslant \lambda_m(x)$. We have
\begin{equation*}
|S_{f,g}(x)|=\tilde{f}_0(x)=\lambda_1(x)\ldots\lambda_m(x)\neq 0
\end{equation*}
and its sign does not change in the interval $I_i$ due to the continuity of $\tilde{f}_0(x)$. Hence each $\lambda_i(x)$ is nonzero and does not change its sign
in $I_i$ since it depends continuously on $x$. It follows that the number of negative eigenvalues of $S_{f,g}(x)$ is constant in $I_i$.
\end{proof}
\begin{proposition}
For $c\in\mathbb R$ with $\tilde{f}_0(c)=0$, there exists an $\varepsilon>0$ such that both $f(x)$ and $g(x)$ are nonzero in
$(c-\varepsilon, c)\cup (c, c+\varepsilon)$. Moreover,
\begin{equation}~\label{diff}
q(S_{f,g}(x))= \left\{ \begin{array}{ll}
q(S_{f,g}(c)) & \textrm{if $f(x)g(x)>0$}\\
q(S_{f,g}(c))+1 & \textrm{if  $f(x)g(x)<0$}
\end{array} \right.
\end{equation}
for $x\in(c-\varepsilon, c)\cup (c, c+\varepsilon)$.
\end{proposition}
\begin{proof}
The first sentence follows from the continuity of $f(x)$ and $g(x)$. We also mention that $\tilde{f}_0(x)$ is nonzero in $(c-\varepsilon, c)\cup (c, c+\varepsilon)$ while $\tilde{f}_1(x)$ is nonzero in $(c-\varepsilon, c+\varepsilon)$, since $\tilde{f}_1(c)\neq 0$.

For simplicity, we denote by $T_{f,g}(x)$ the principal submatrix $S_{f,g}(x)_{\{2,3,\ldots, m\}}$ of $S_{f,g}(x)$. By Proposition~\ref{minor},
\begin{equation*}
|T_{f,g}(x)|=D_{m-1}(x)=\tilde{f}_1(x).
\end{equation*}
Since the sign of $\tilde{f}_1(x)$ does not change in $(c-\varepsilon, c+\varepsilon)$, using a similar argument as in Proposition~\ref{constant},
we can prove that $q(T_{f,g}(x))$ is constant in $(c-\varepsilon, c+\varepsilon)$. It follows from Lemma~\ref{lem1} that for $x\in(c-\varepsilon, c)\cup (c, c+\varepsilon)$,
\begin{displaymath}
q(S_{f,g}(x))= \left\{ \begin{array}{ll}
q(T_{f,g}(x)) & \textrm{if $\tilde{f}_0(x)\tilde{f}_1(x)>0$}\\
q(T_{f,g}(x))+1 & \textrm{if  $\tilde{f}_0(x)\tilde{f}_1(x)<0$}
\end{array} \right.
\end{displaymath}
and
\begin{equation*}
q(S_{f,g}(c)=q(T_{f,g}(c))=q(T_{f,g}(x)).
\end{equation*}
Thus (\ref{diff}) is an immediate result since $f(x)g(x)>0$ is equivalent to $\tilde{f}_0(x)\tilde{f}_1(x)>0$ for $x\in(c-\varepsilon, c)\cup (c, c+\varepsilon)$.
\end{proof}

Set $q(S_{f,g}(x))=q_i$ in $I_i$ and $q(S_{f,g}(a_j))=t_j$ with $0\leqslant i\leqslant k$ and $1\leqslant j\leqslant k$. The following corollary is implied by above proposition.
\begin{corollary}~\label{cor} For $1\leqslant i\leqslant k$,
\begin{equation*}
\begin{split}
t_i= &\left\{ \begin{array}{ll}
q_{i-1} & \textrm{if $f(x)g(x)>0$ in some interval $(a_i-\varepsilon, a_i)$}\\
q_{i-1}+1 & \textrm{if $f(x)g(x)<0$ in some interval $(a_i-\varepsilon, a_i)$}
\end{array} \right.\\
=&\left\{ \begin{array}{ll}
q_{i} & \textrm{if $f(x)g(x)>0$ in some interval $(a_i, a_i+\varepsilon)$}\\
q_{i}+1 & \textrm{if $f(x)g(x)<0$ in some interval $(a_i, a_i+\varepsilon)$}
\end{array} \right.
\end{split}
\end{equation*}
\end{corollary}
If $g(x)$ is the derivative of $f(x)$, i.e. $g(x)=f^{\prime}(x)$, $\tilde{f}_0(x)$ has exactly the same roots as $f(x)$ with multiplicity 1. As a real polynomial function, $f^2(x)$ reaches its minimum at $a_i$.
It is a simple fact that $(f^{2}(x))^{\prime}=2f(x)f^{\prime}(x)>0$ in some interval $(a_i, a_i+\varepsilon)$ and $2f(x)f^{\prime}(x)<0$ in some interval $(a_i-\eta, a_i)$. Hence we can prove Sturm's theorem by Corollary~\ref{cor}.
\begin{theorem}[Sturm's Theorem]
Let $a_1<a_2<\dotsb<a_k$ be the real roots of $f(x)\in\mathbb R[x]$, and $g(x)$ be the derivative of $f(x)$.
\begin{itemize}[align=left]
\item[(i)] $q(S_{f,g}(x))$ is constant on each interval $(a_{i-1}, a_i]$ with value $q_i$ for $i=1,2,\ldots,k$, where $a_0=-\infty$.
\item[(ii)] $q_i-q_{i+1}=1$ for $0\leqslant i\leqslant k-1$.
\item[(iii)] For any $a<b$, the number of distinct real roots of $f(x)$ in the interval $(a, b]$ is $q(S_{f,g}(a))-q(S_{f,g}(b))$.
\end{itemize}
\end{theorem}
\begin{remark}
One can see that the above theorem covers the original one by Theorem~\ref{var} while the condition that $a, b$ are not multiple roots of $f(x)$ is unnecessary in this modified version of Sturm's Theorem.
\end{remark}

\end{document}